\documentclass{article}
\usepackage{amsmath,amsthm,amssymb,mathrsfs}
\usepackage[utf8]{inputenc}\usepackage[pdftex]{graphicx}\usepackage{hyperref}
\usepackage[headsep=0.8in, textheight=700pt]{geometry}
\usepackage{parskip}
\usepackage{float}
\usepackage{listings}
\usepackage{eurosym}
\usepackage[font=small,labelfont=bf]{caption}
\usepackage{subcaption}
\newcommand{\xib}{\hat{\xi}_k^{\mathcal{B}}}
\newcommand{\db}{\hat{\delta}_k^{\mathcal{B}}}
\theoremstyle{plain}

\newlength{\figwidth} 

\hypersetup{pdfinfo={
	Title={Reducing MSE in estimation of heavy tails: a Bayesian approach},
	Subject={Bayesian bias-reduced estimation of the extreme value index for heavy-tailed distributions},
	Author={Gaonyalelwe Maribe, Andréhette Verster and Jan Beirlant},
	Keywords={Extended Pareto Distribution, extreme value index, tail estimation, Bayesian parameter estimation, bias reduction, posterior simulation} } }
\DeclareRobustCommand{\officialeuro}{%
	\ifmmode\expandafter\text\fi
	{\fontencoding{U}\fontfamily{eurosym}\selectfont e}}

\title{Reducing bias and MSE in estimation of heavy tails:\\ a Bayesian approach}
\author{Gaonyalelwe Maribe\thanks{University of the Free State. This work is based on the research supported wholly/in part by the National Research Foundation of South Africa (Grant Number 102628). The Grantholder acknowledges that opinions, findings and conclusions or recommendations expressed in any publication generated by the NRF supported research is that of the author(s), and that the NRF accepts no liability whatsoever in this regard.} \and Andréhette Verster\thanks{University of the Free State} \and Jan Beirlant\thanks{KU Leuven and University of the Free State}}

\pdfoutput=1
\begin{document}
	\maketitle	
\begin{abstract}
Bias reduction in tail estimation has received considerable interest in extreme value analysis.  Estimation methods that minimize the bias while keeping the mean squared error (MSE) under control, are especially useful when applying classical methods such as the Hill (1975) estimator. In Caeiro {\it et al.} (2005) minimum variance reduced bias estimators of the Pareto tail index were first proposed where the bias is reduced without increasing the variance with respect to the Hill estimator. This method is based on adequate external estimation of a pair of second-order parameters. Here we revisit this problem from a Bayesian point of view starting from the extended Pareto distribution (EPD) approximation to excesses over a high threshold, as developed in Beirlant {\it et al.} (2009) using maximum likelihood (ML) estimation. Using asymptotic considerations, we derive an appropriate  choice of priors leading to a Bayes estimator for which the  MSE curve  is a weighted average  of the Hill and EPD-ML MSE curves for a large range of thresholds, under the same conditions as in Beirlant {\it et al.} (2009).  A  similar result is obtained for tail probability estimation. Simulations show surprisingly good MSE performance with respect to the existing estimators.
		\end{abstract}

	
\section{Introduction}
In this paper we consider the estimation of the extreme value index $\xi$ and tail probabilities $P(X>x)$ for $x$ large, on the basis of  independent and identically distributed observations $X_1, X_2, \ldots, X_n$  which follow a  Pareto-type distribution with right tail function (RTF) given by
\begin{equation}
\bar{F}(x) = 1-F(x) = P(X>x) = x^{-1/\xi} \ell (x)
\label{Patype}
\end{equation}
where $\ell$ is a slowly varying function at infinity, {\it i.e.}
$$
{\ell (ty) \over \ell (t)} \to 1, \mbox{ as } t \to \infty, \mbox{ for every  } y>1.
$$
The most famous estimator of $\xi$ was first derived by Hill (1975) as a maximum likelihood (ML) estimator using the approximation 
\begin{equation}
\bar{F}(ty)/\bar{F}(t) \approx y^{-1/\xi}, t \mbox{ large}
\label{Pot}
\end{equation}
to the RTF  of the excesses $X/t |X>t$ over a large threshold $t$  by a simple Pareto distribution with RTF $y^{-1/\xi}$, and setting $t= X_{n-k,n}$ where $X_{1,n}\leq X_{2,n} \leq \ldots \leq X_{n,n}$:
\begin{equation}
H_{k,n} = {1 \over k}\sum_{j=1}^k \log {X_{n-j+1,n} \over X_{n-k,n}}.
\label{Hill}
\end{equation}
A simple estimator of a tail probability $P(X>x)$ with $x$ large, introduced in Weissman (1978),\\ is then obtained from \eqref{Pot} setting $ty = x$ and estimating $P(X>t)$ by the empirical proportion $k/n$: 
\begin{equation}
\hat{p}_{x,k} = {k \over n}\left( {x \over X_{n-k,n}}\right)^{-1/H_{k,n}}.
\label{Weiss}
\end{equation}
In practice, a way to verify the validity of model \eqref{Patype} is to check whether the Hill estimates are stable as a function of $k$. However in most cases the stability is not visible, which can be explained by slow convergence in \eqref{Pot}. For this reason bias reduced estimators have been proposed which lead to plots that are much more horizontal in $k$ which facilitates the analysis of a practical case to a great extent.   Here we can refer to Peng (1998), Beirlant {\it et al.} (1999, 2008), Feuerverger and Hall (1999), Caeiro {\it et al.} (2005, 2009) and Gomes {\it et al.} (2007) for bias-reduced estimators based on functions of the top $k$ order statistics. Several of these methods focus on the distribution of log-spacings of high order statistics.\\\\ Beirlant {\it et al.} (2009) proposed to build a more flexible model capable of capturing the deviation between the true excess RTF  $\bar{F}(ty)/\bar{F}(t)$ and the asymptotic Pareto model. For a heavy tailed distribution \eqref{Patype}, this deviation can be parametrized using a power series expansion (Hall, 1982), or more generally via second-order slow variation (Bingham {\it et al.}, 1987). 
More specifically in Beirlant {\it et al.} (2009) the subclass $\mathcal{F}(\xi,\tau)$ of the Pareto-type tails \eqref{Patype} was considered satisfying
\begin{eqnarray}
\bar{F}(x) = C x^{-1/\xi}\left( 1+\xi^{-1}\delta (x)\right),
\label{Sclass}
\end{eqnarray}
with $\delta (x)$ eventually nonzero and of constant sign such that $|\delta(x)|=x^\tau \ell_{\delta}(x)$ with $\tau <0$ and $\ell_{\delta}$ slowly varying.  
It was shown that under $\mathcal{F}(\xi,\tau)$ one has as $t \to \infty$
$$
\sup_{y \geq 1}\left|{\bar{F}(ty)\over \bar{F}(t)} - \bar{G}_{\xi,\delta,\tau}(y) \right|
= o\left(|\delta (t)| \right)
$$
with $\bar{G}_{\xi,\delta,\tau}$ the RTF of the extended Pareto distribution (EPD) 
\begin{equation}
\bar{G}_{\xi,\delta,\tau}(y)=\{y(1+\delta-\delta y^{\tau})\}^{-1/\xi}, \ \ \ y>1,
\label{EPD}
\end{equation}
with $\tau <0 < \xi$ and $\delta > \max (-1,1/\tau)$.
This shows that the EPD improves the approximation \eqref{Pot} with an order of magnitude.
Then maximum likelihood (ML) estimation of the parameters ($\xi,\delta$) based on a set of  excesses 
$\left(Y_{j,k}:=X_{n-j+1,n}/X_{n-k,n},\; j=1,\ldots,k \right)$ was used to obtain a bias reduced estimator of $\xi$. Bias reduction of the Weissman estimator of tail probabilities can analogously be obtained using 
\begin{equation}
\hat{p}^{EP}_{x,k} = {k \over n}\bar{G}_{\hat\xi _k,\hat\delta _k,\hat\tau}\left( {x \over X_{n-k,n}}\right),
\label{PEP}
\end{equation}
where $(\hat\xi _k,\hat\delta _k)$ denote the ML estimators based on the EPD model, and where $\hat\tau$ is a consistent estimator of $\tau$, to be specified below, which was shown not to affect the asymptotic distribution of $(\xi,\delta)$.

\vspace{0.3cm}\noindent
Here we investigate the possibilities of using the Bayesian methodology when modelling the distribution of the vector of excesses ${\bf Y}_k := \left( Y_{j,k}, j=1,\ldots,k \right)$
with an EPD. In section 2 we show that a normal prior on $\delta$ with zero mean and  variance $\sigma^2_{k,n}$, depending in an appropriate way on $k$ and $n$, leads to interesting MSE results of the posterior mode estimators for $\xi$ and of the corresponding estimators of $P(X>x)$ following \eqref{PEP}. In section 3 we discuss the implementation of the Bayesian estimates. In section 4 we consider the finite sample behaviour of this Bayesian approach and consider some practical cases. 

\section{Bayesian estimation of the EPD parameters}

ML estimation of the EPD parameters ($\xi,\delta$), given a value of $\tau$, follows by maximizing the log-likelihood
\begin{equation}
{1 \over k}l(\xi,\delta |{\bf y}) =  -\log \xi
-(\frac{1}{\xi}+1){1 \over k}\sum_{j=1}^k\left[\log y_{j,k} + \log(1+\delta\{1-y_{j,k}^{\tau}\})\right]
+{1 \over k}\sum_{j=1}^k\log\left(1+\delta\{1-(1+\tau)y_{j,k}^{\tau}\}\right).
\label{logl}
\end{equation}
In the Bayesian framework, the log-posterior can be written as
 \begin{equation}
 {1 \over k} \log \pi (\xi,\delta|{\bf y}) = 
{1 \over k}l(\xi,\delta |{\bf y}) + {1 \over k}\log \pi (\xi,\delta),
\label{logpost}
\end{equation}
where $\pi (\xi,\delta)$ denotes the prior density. Here   we assign a maximal data information  (MDI) prior to $\xi$, which for a general parameter $\theta$ is defined as
$
\pi(\theta)\propto\exp(E(\log f({\bf Y}|\theta)))$.
Beirlant {\it et al.} (2004) derived that the MDI for a Pareto distribution  is given by
\begin{equation}\pi(\xi)\propto \frac{e^{-\xi}}{\xi}.
\label{pixi}
\end{equation}
\\
Next, the prior on $\delta$ is taken to be a  normal distribution with mean 0 and variance $\sigma^2_{k,n}$, depending on  
$k$, and left truncated  in order to comply with the restriction
 $\delta > \max (-1,1/\tau)$:
\begin{equation}
\pi (\delta) = {1 \over \sqrt{2\pi}\sigma_{k,n}}e^{-{1 \over 2}{\delta^2 \over \sigma^2_{k,n}}}/ \left( 1-\Phi (\max (-1,\tau^{-1})/\sigma_{k,n})\right).
\label{pidelta}
\end{equation} 

\vspace{0.4cm}\noindent
If $F$ satisfies $\mathcal{F}(\xi,\tau)$, it is shown in  Beirlant {\it et al.} (2009) that $Q(1-x^{-1})$ ($x>1$), with 
$Q(p) = \inf \{x: F(x) \geq p \}$ ($p \in (0,1)$), satisfies
\begin{equation}
Q(1-x^{-1}) = C^{\xi} x^{\xi}\left(1+ a(x) \right)
\label{Q}
\end{equation}
with $a(x) = \delta (Q(1-x^{-1}))\{1+o(1)\}= \delta (C^{\xi} x^{\xi})\{1+o(1)\}$ as $x \to \infty$. In particular
$a$ is eventually nonzero and of constant sign and $|a(x)|= x^{\rho}\ell_{a}(x)$ with $\ell_a$ slowly varying and $\rho=\xi\tau$. 

Then it was shown that if $F$ satisfies $\mathcal{F}(\xi,\tau)$, and $\sqrt{k} a(n/k) \to \lambda \in \mathbb{R}$ and $\hat{\rho}_n = \rho + o_p(1)$ as $k,n \to \infty$ and $k/n \to 0$, the following asymptotic results hold for the EPD-ML estimator $\hat{\xi}^{ML}_{k,n}$ and $H_{k,n}$:
\begin{eqnarray}
\sqrt{k} \left( \hat{\xi}^{ML}_{k,n}-\xi \right) & \to_d & 
\mathcal{N}\left(0, \xi^2 \left({1-\rho \over \rho} \right)^2
 \right),
 \label{limML} \\
 \sqrt{k} \left( H_{k,n}-\xi \right) & \to_d & 
\mathcal{N}\left(\lambda {\rho \over 1-\rho} , \xi^2
 \right).
  \label{limH} 
\end{eqnarray}
From these results it follows that for the smallest values of $k$ ({\it i.e.} $\sqrt{k} a(n/k) \to 0$) the Hill estimator is asymptotically unbiased, while for increasing values of $k$ ({\it i.e.} $\sqrt{k} a(n/k) \to \lambda \neq 0$) it is biased. In this region the bias reduced estimator $ \hat{\xi}^{ML}_{k,n}$  still has asymptotic bias 0, but its variance is increased by a factor $((1-\rho)/\rho)^2$ compared to $H_{k,n}$.    
\\\\
In the Appendix we derive that the first order approximations ($\xib,\db$) of the Bayesian estimators are given by
  \begin{eqnarray*}
  \xib & =& H_{k,n} +\db \left( 1-E_{k,n}(\tau) \right) ,\\
    \db &=& \frac{1-H_{k,n}\tau}{D^{\mathcal{B}}_{k,n}} \left(E_{k,n}(\tau) - {1 \over H_{k,n}\tau} \right)
  \end{eqnarray*}
  where  
 $$
 E_{k,n}(s) = {1 \over k}\sum_{j=1}^k Y_{j,k}^{s},\; s < 0,
  $$
  and 
 $$
  D^{\mathcal{B}}_{k,n} =
  {\xib \over k\sigma^2_{k,n}} -
  \left( 
  1- 2(1-\xib \tau) E_{k,n}(\tau)+ (1-2\xib\tau-\xib\tau^2)E_{k,n}(2\tau) -\tau (1- E_{k,n}(\tau))E_{k,n}(\tau)
  \right).
  $$
 These expressions are identical to the asymptotic EPD-ML estimators derived in Beirlant {\it et al.} (2009) except for the extra term ${\xib\over k\sigma^2_{k,n}}$ in the expression of  $D^{\mathcal{B}}_{k,n}$. 
   As an external estimator of $\tau$ we use $\hat{\tau} = \hat{\rho}/H_{k,n}$ where $\hat{\rho}$ is a consistent estimator of $\rho$ (see for instance $\rho_{Fraga}$ from Fraga Alves {\it et al.}, 2003). The following result is derived in the Appendix. 
  \\\\
{\bf Theorem.} {\it Let $F \in \mathcal{F}(\xi,\tau)$ and assume $k\sigma_{k,n}^2 \to \mu >0$ as $k,n \to \infty$, $k/n \to 0$ and $\sqrt{k} a(n/k) \to \lambda$. Then $\Xi_{k,n}:=\sqrt{k} \left( \xib-\xi \right)$ is asymptotically normal  with asymptotic mean and variance given by
  \begin{eqnarray}
 E_{\infty} (\Xi_{k,n}) &=& {\lambda \rho \over 1-\rho}{\zeta \over \zeta + \rho^4}, 
 \label{Abias}\\
  Var_{\infty} (\Xi_{k,n})&=& {\xi^2 \over (1+ \zeta \rho^{-4})^2}\left( \left({1-\rho \over \rho} \right)^2
  + {\zeta^2 \over \rho^8}+ 2{\zeta \over \rho^4}\right).
  \label{Avar}
 \end{eqnarray}
 }
 
 \vspace{0.3cm}\noindent
 Minimizing $MSE_{\infty}(\Xi_{k,n})=  E^2_{\infty} (\Xi_{k,n})
 +  Var_{\infty} (\Xi_{k,n})$ with respect to $\zeta$ leads to
\begin{equation}
 \zeta_{opt}= {\xi^2(1-2\rho) \over \lambda^2},
 \label{zetaopt}
\end{equation}
 from which, with $\zeta_{opt}= \xi^2 (1-2\rho)(1-\rho)^2/\mu_{opt}$,
\begin{equation}
 \mu_{opt}= (1-\rho)^2
\lambda^2. 
\label{muopt}
\end{equation}
Hence, we are lead to choosing 
\begin{equation}
\sigma^2_{opt,k,n} = (1-\rho)^2 a^2(n/k)
\label{sigmaopt}
\end{equation}
for the prior variance since $\mu = \lim_{k,n \to \infty} k\sigma^2_{k,n}$ and $\lambda^2 = \lim_{k,n \to \infty} ka^2(n/k)$. \\\\
Note that using \eqref{zetaopt}, one obtains from \eqref{Abias} and \eqref{Avar} that 
\begin{eqnarray*}
E^{opt}_{\infty} (\Xi_{k,n}) &=& {\lambda \rho\over 1-\rho}
{\xi^2 (1-2\rho) \over \xi^2 (1-2\rho)+ \lambda^2 \rho^4 },\\
Var^{opt}_{\infty} (\Xi_{k,n}) &=& {\xi^2 \lambda^4 \over (\lambda^2\rho^4+ \xi^2(1-2\rho))^2}
\left\{(1-\rho)^2\rho^6 +  \lambda^{-4}\xi^4(1-2\rho)^2 + 2\lambda^{-2}\xi^2(1-2\rho)\rho^4\right\},
\end{eqnarray*}
from which 
\begin{eqnarray}
MSE^{opt}_{\infty}(\Xi_{k,n}) &=&
\frac{\xi^4(1-2\rho)^2[\xi^2 + \lambda^2{\rho^2 \over (1-\rho)^2}]+2 \xi^2\rho^4\lambda^2(1-2\rho)[\xi^2]+\lambda^4\rho^8 [\xi^2{(1-\rho)^2\over \rho^2}] }{\xi^4(1-2\rho)^2 + 2 \xi^2\rho^4\lambda^2(1-2\rho)+ 
\lambda^4\rho^8} \label{wavg} \\
&=& \xi^2 + \xi^2\rho^2{(1-2\rho)\over (1-\rho)^2}\, \lambda^2 \, \frac{\left(\xi^2(1-2\rho)+ \lambda^2\rho^4(1-\rho)^2 \right)}{(\xi^2(1-2\rho)+\lambda^2\rho^4)^2}. \label{xisqdiff}
\end{eqnarray}
From \eqref{wavg} it follows that the asymptotic MSE of this optimal Bayesian estimator is a weighted average of the asymptotic MSEs of the Hill and EPD-ML estimators. Moreover, since the right hand side of  \eqref{xisqdiff} is an increasing function in $\lambda^2$ it follows that
$$
MSE^{opt}_{\infty}(\Xi_{k,n})  \leq \lim_{\lambda \to \infty}MSE^{opt}_{\infty}(\Xi_{k,n}) = MSE_{\infty}\left(\sqrt{k}(\hat{\xi}_{k,n}^{ML}-\xi)\right).
$$ 
Also, expanding the the right hand side of  \eqref{xisqdiff} for $\lambda^2 \to 0$ leads to 
$$
MSE^{opt}_{\infty}(\Xi_{k,n}) = \xi^2 + \lambda^2 {\rho^2 \over (1-\rho)^2}(1+o(1)).
$$
We can conclude that the asymptotic MSE of the optimal Bayes estimator is uniformly smaller than the MSE of the EPD-ML estimator as given in \eqref{limML}, while for smaller $\lambda$ this asymptotic MSE follows the asymptotic MSE of the Hill estimator, given in \eqref{limH}, up to terms of order $\lambda^2$. Hence with the choice \eqref{sigmaopt} the Bayesian estimator automatically follows the better of the two existing estimators as a function of $\lambda$ or $k$.
\\\\
Replacing ($\hat\xi _k,\hat\delta _k$) by ($\hat\xi ^{\mathcal{B}}_k,\hat\delta ^{\mathcal{B}}_k$) in  $\hat{p}^{EP}_{x,k}$, it follows from the proof of Theorem 5.2 in Beirlant {\it et al.} (2009) that the resulting tail probability estimator $\hat{p}^{\mathcal{B}}_{x,k}$ satisfies the following asymptotic result under the conditions of the Theorem: \\
when $p_n = P(X>x_n)$ satisfies $np_n/k \to 0$ and $\log(np_n)/\sqrt{k} \to 0$, then 
$${\sqrt{k} \over \log (k/(np_n))}\xi ( {\hat{p}^{\mathcal{B}}_{x_n,k} \over p_n}-1)$$ is asymptotically normal with the same limit distribution as in the Theorem. Hence the  asymptotic MSE behaviour  for the tail probability estimator has the same characteristics as  the tail index estimator.
\\\\
In the popular special case which contains most Pareto-type distributions (Hall, 1982) where for some constant $C$
\begin{equation}
x^{-\rho}a(x) \to_{x \to \infty} C,
\label{Hall}
\end{equation} 
the limit $C$ can be incorporated in  $\lambda$ so that under \eqref{Hall}  the variance of the normal prior on $\delta$ can be taken as
\begin{equation}
\sigma^2_{*,k,n} = \left({k \over n}\right)^{-2\hat\rho}.
\label{sigmapract}
\end{equation}
This is the choice of $\sigma^2_{k,n}$ we will use throughout in practice, with $\rho$ estimated by the method proposed in Fraga Alves (2003). We will also study the sensitivity of the method when replacing $\rho$ by $-1$, as used for instance in Beirlant {\it et al.} (2002).
\\\\
Denoting the Bayesian estimators of ($\xi,\delta$) using \eqref{sigmapract} by $(\hat{\xi}_{k,n}^*,\hat{\delta}_{k,n}^*)$, as in \eqref{PEP} we then find the corresponding estimator for tail probabilities $P(X>x)$:
\begin{equation}
\hat{p}^{*}_{x,k} = {k \over n}\bar{G}_{(\hat{\xi}^* _{k,n},\hat\delta^* _{k,n},\hat\rho/H_{k,n})}\left( {x \over X_{n-k,n}}\right).
\label{PBayes}
\end{equation}

\section{Implementation of the Bayes estimators}
Bayesian inference has a great advantage of incorporating in a unified way, any meaningful piece of information in describing our model parameters. Here we use the objective mathematical information that with decreasing $k$ (or increasing threshold $t$) the $\delta$ parameter becomes smaller, expressed by the variance of the normal prior on this parameter which is taken to be of order $(k/n)^{-2\rho}$. We rely on modern Markov Chain Monte Carlo (MCMC) methods to assist  in approximating posterior distributions. In order to make inference on $(\xi,\delta)$ we take the mode of the posterior samples, and the accuracy of this inference is described by the posterior distribution itself through the highest posterior density (HPD) region. The HPD for $\xi$ represents a set of most probable values of the $\hat{\xi}$, constituting $100(1-\alpha)\%$ of the posterior mass. The HPD also has the characteristic that the density within the HPD region is never lower than the values outside.

The Bayesian estimations are implemented via OpenBUGS (an open source version of BUGS) in the \textbf{R} (R Core Team, 2014) statistical software.
\\\\
In OpenBUGS we  need to specify the EPD model likelihood, the parameter priors and starting values, then a Markov chain simulation is automatically implemented for the resulting posterior distribution given in (\ref{logpost}). The EPD is however not included in the standard distributions list available in OpenBUGS, and we therefore indirectly implement the EPD model likelihood using a Poisson distribution.\\ \\
Let $l_i=\frac{1}{k}l(\xi,\delta|y_i)$, where $\frac{1}{k}l(\xi,\delta|y_1,...,y_n)$ is the EPD log-likelihood defined in (\ref{logl}). The EPD model likelihood can be written as\begin{equation}
f(\textbf{y}|\xi,\delta)=\prod_{i=1}^{n}e^{l_i}=\prod_{i=1}^{n}\frac{e^{-(-l_i)}(-l_i)^0}{0!}=\prod_{i=1}^{n}f_{P}(0;-l_i).
\end{equation}
The resulting density $f_{P}(0;-l_i)$  therefore is a product of Poisson distributed pseudo random variables with mean equal to the EPD log-likelihood and with all observed values set equal to zero. This is known as the `Zero trick' (Lunn {\it et al.}, 2012). A positive constant $C$ is added to the mean to ensure that the mean of the pseudo random variables is positive. The resulting likelihood becomes\begin{equation}
f(\textbf{y}|\xi,\delta)=\prod_{i=1}^{n}\frac{e^{-(-l_i+C)}(-l_i+C)^0}{0!}=\prod_{i=1}^{n}f_{P}(0;-l_i+C);
\end{equation}where $C$ is chosen such that $-l_i+C>0 \mbox{ for } i=1,2,...,n$.\\ \\
In BUGS the full probability model needs to be defined and hence all prior distributions need to be proper ({\it i.e.} integrate to 1). The left truncated normal prior on $\delta$ is proper and can be specified directly in OpenBUGS. To ensure the MDI Pareto prior on $\xi$ is proper, we write it as
\begin{equation}\pi(\xi)\propto \frac{e^{-\xi}}{\xi} \propto \frac{1}{\Gamma(\alpha)}\xi^{\alpha-1}e^{-\xi}. 
\label{pixi}
\end{equation}and specify it in OpenBUGS using a Gamma distribution with a scale parameter equal to 1 and a shape parameter equal to 0.0001.
\\\\
Values of $\rho$ close to 0 have to be avoided, and so we put a restriction on $\rho_{Fraga}$ in our implementation, using $\hat{\rho}= \min (-0.5, \rho_{Fraga})$. Finally we also smooth the Bayesian estimates as a function of $k$ with a moving average with a window width of 5. 

\section{Simulations and practical case studies}
We performed a simulation study, taking 1000 repetitions of samples of size $n=500$ studying the finite sample behaviour of $\hat{\xi}^*_{k,n}$ and $\hat{p}^{*}_{x,k}$ with $x=Q(1-1/500)$  for different distributions, using the estimator $\hat\rho$ proposed by Fraga Alves {\it et al.} (2003) and the results are shown in figures 1 to 3. The bias, variance and MSE are plotted as a function of $k$. In case of the tail probability estimators we consider the relative error in the variance and MSE.\\
\\
We compute the HPD using a direct approach of taking the shortest probability interval, given a $1-\alpha$ coverage based on the simulations. We do this by ordering the $m$ simulation draws and then taking the shortest interval that contains $(1-\alpha)*m$ of the draws.\\
The following distributions are used:
\begin{itemize}
 \item
{\it The Fr\'echet distribution} with  $\bar{F}(x)= 1-\exp(-x^{-1/\xi})$ taking $\xi = 0.5$ in which case $\rho=-1$.
\item
{\it The Burr distribution} with  $\bar{F}(x)= (1+x)^{-4/3}$ so that $\xi = 0.75$ and $\rho=-0.75$.
\item
{\it The loggamma distribution} with $\bar{F}(x) \sim constant \times x^{-2}(\log x)^3$ so that $\xi = 0.5$, which does not belong to the class  $\mathcal{F}(\xi,\tau)$.
\end{itemize}
 	
 		\begin{figure}[H]
 			\centering
 			\begin{subfigure}[h]{0.47\textwidth}
 			\includegraphics[width=80mm]{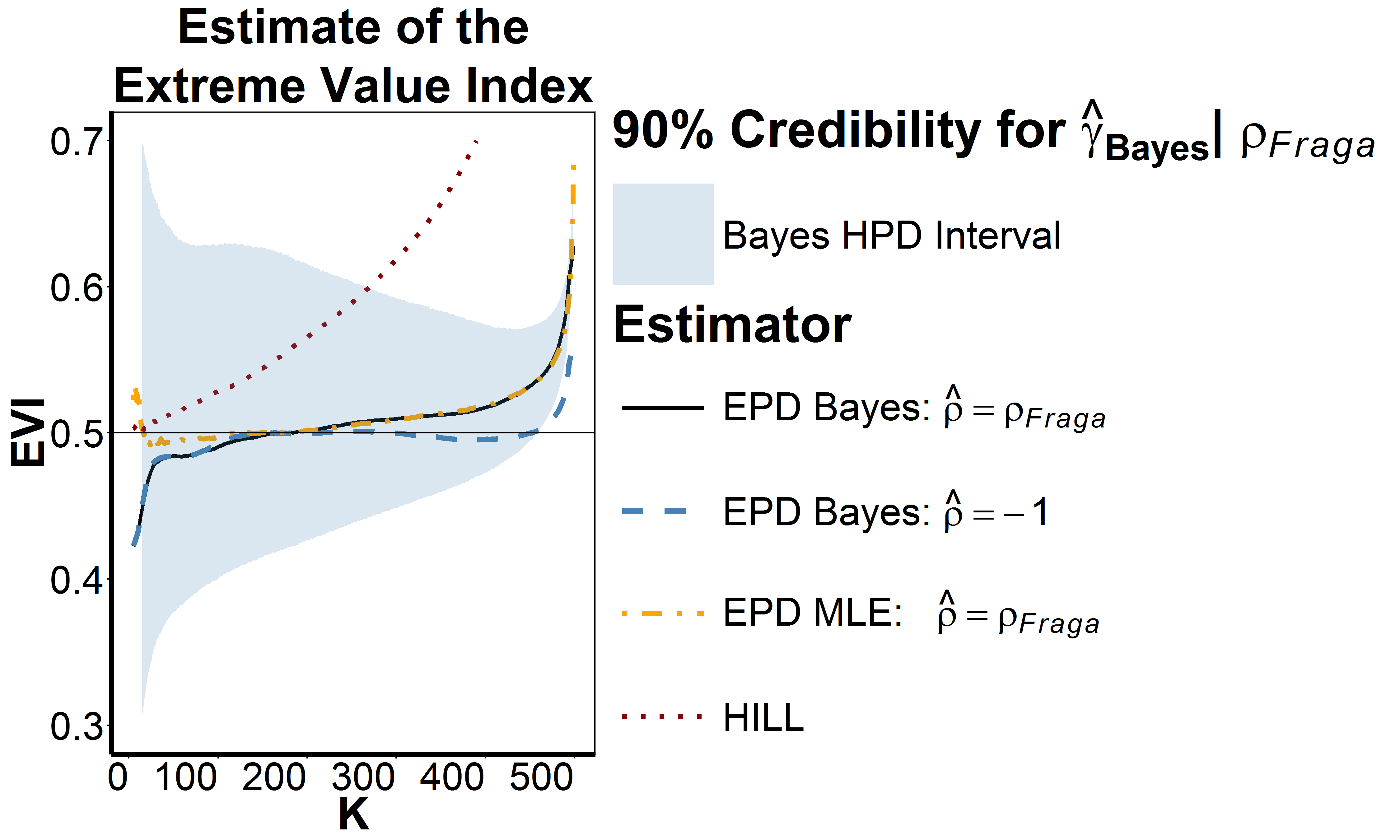}
 			\end{subfigure}
 			\hspace{\fill}
 			\begin{subfigure}[h]{0.45\textwidth}
 			\includegraphics[width=85mm]{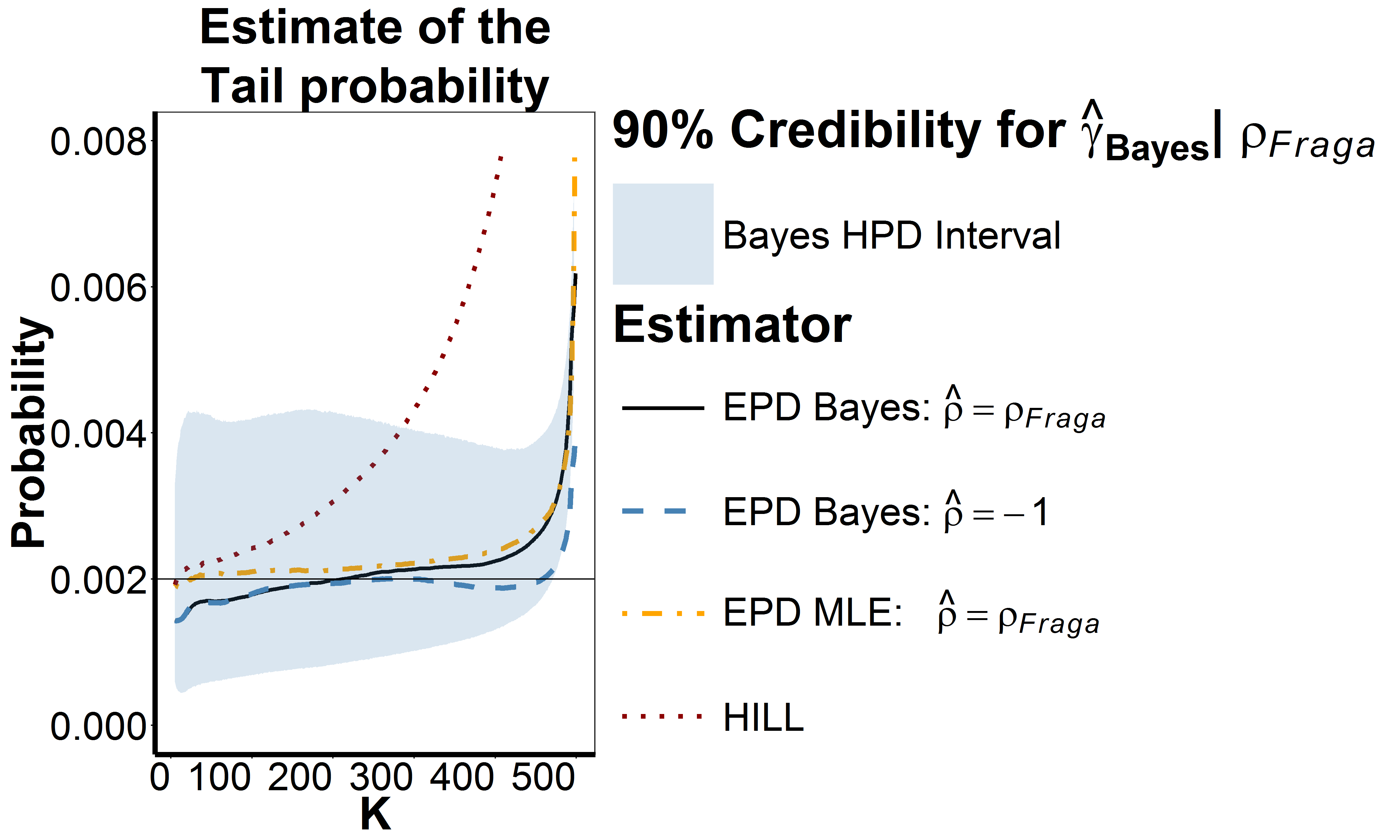}
 			\end{subfigure}
 			\bigskip
 			
 			\begin{subfigure}[h]{0.45\textwidth}
 			\includegraphics[width=75mm]{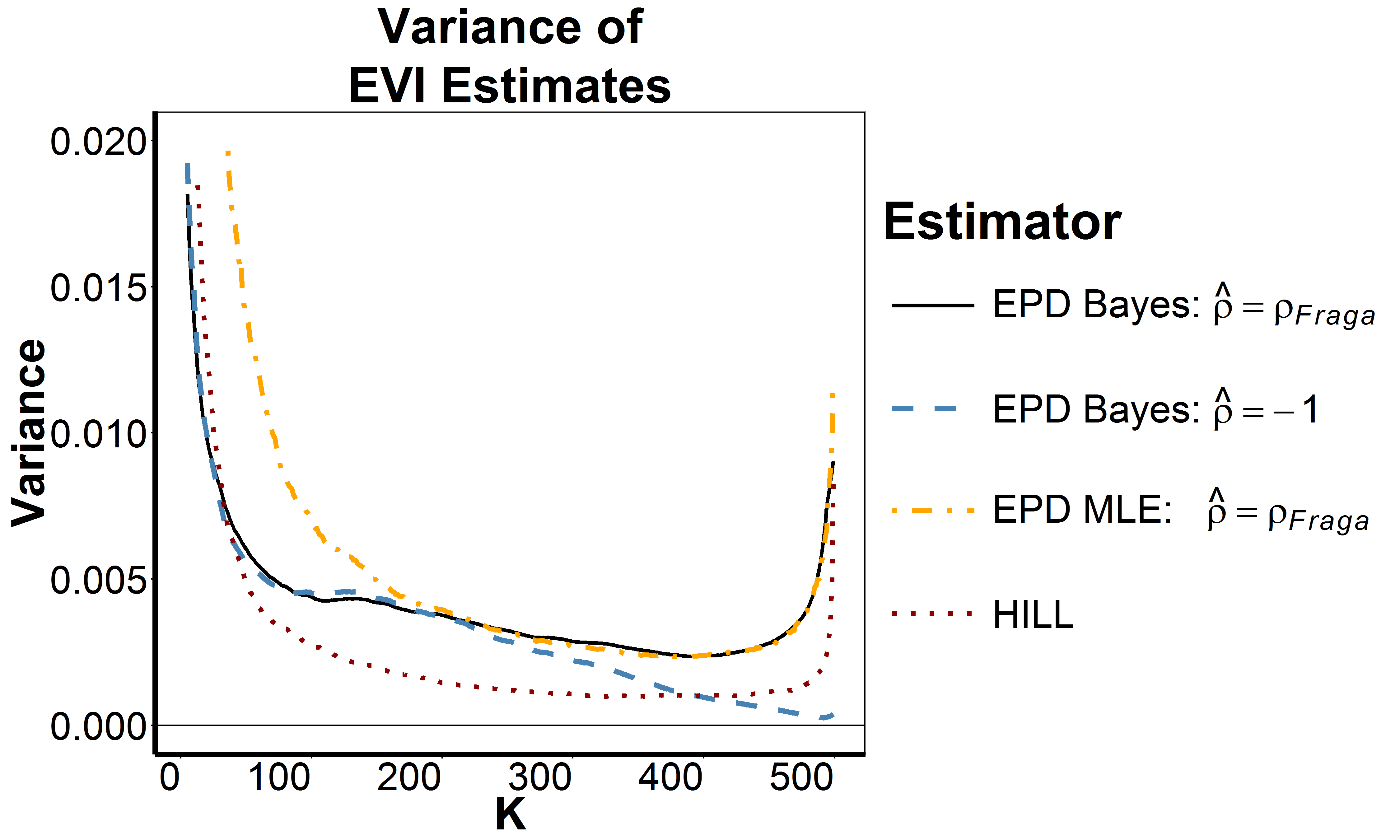}
 			\end{subfigure}	
 			\hspace{\fill}
 			\begin{subfigure}[h]{0.45\textwidth}
 			\includegraphics[width=75mm]{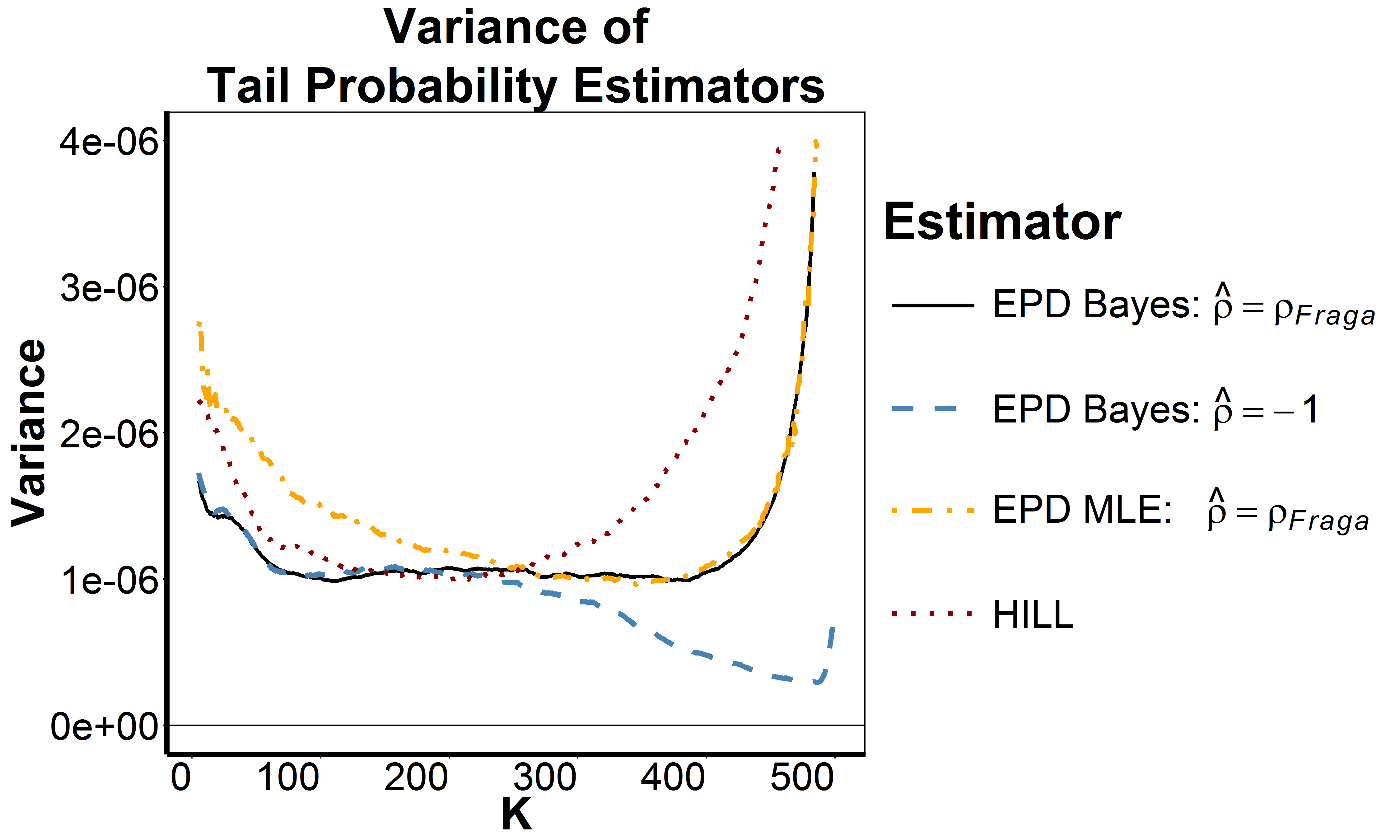}
 			\end{subfigure}	
 			\bigskip
 			
 			\begin{subfigure}[h]{0.45\textwidth}
 			\includegraphics[width=75mm]{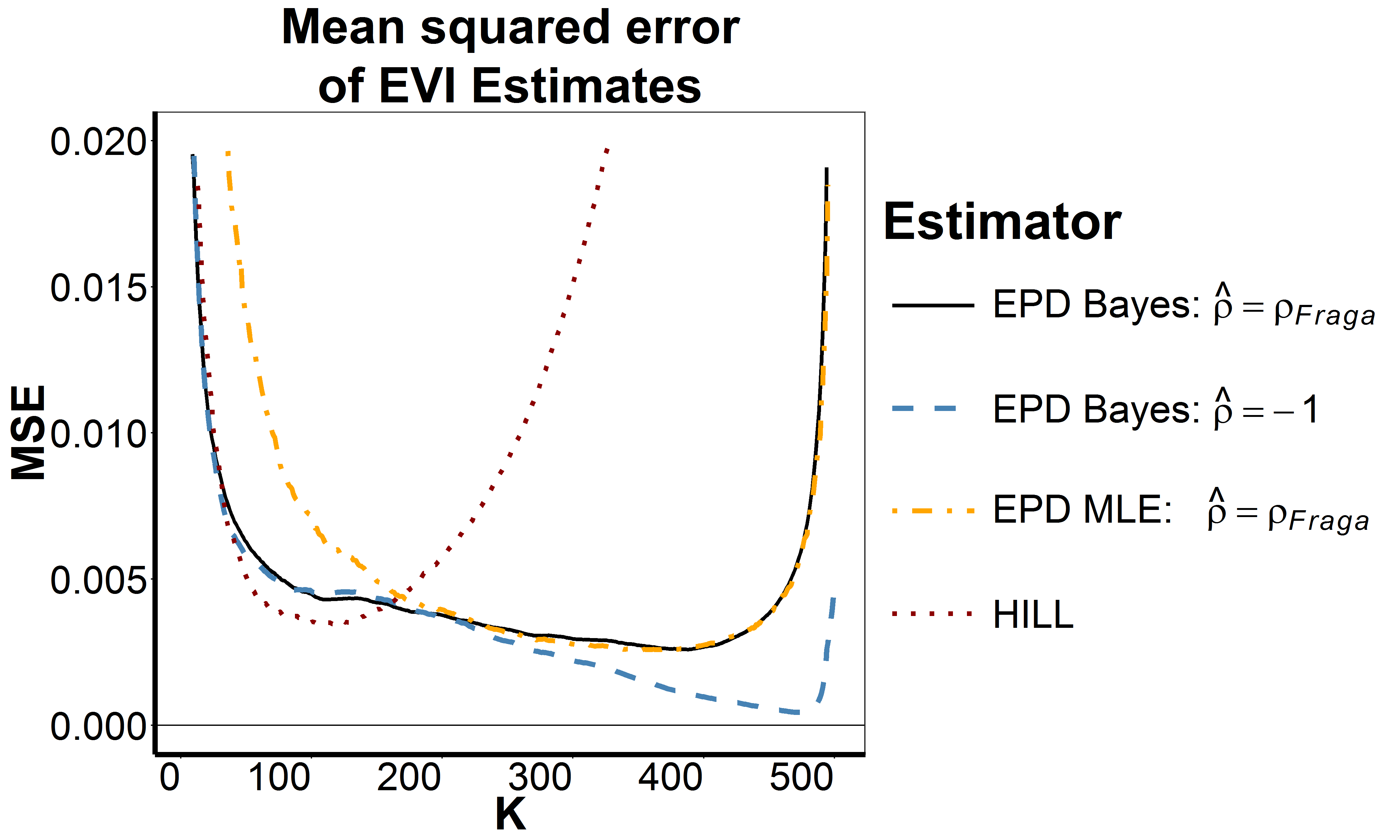}
 			\end{subfigure}	
 			\hspace{\fill}
 			\begin{subfigure}[h]{0.45\textwidth}
 			\includegraphics[width=75mm]{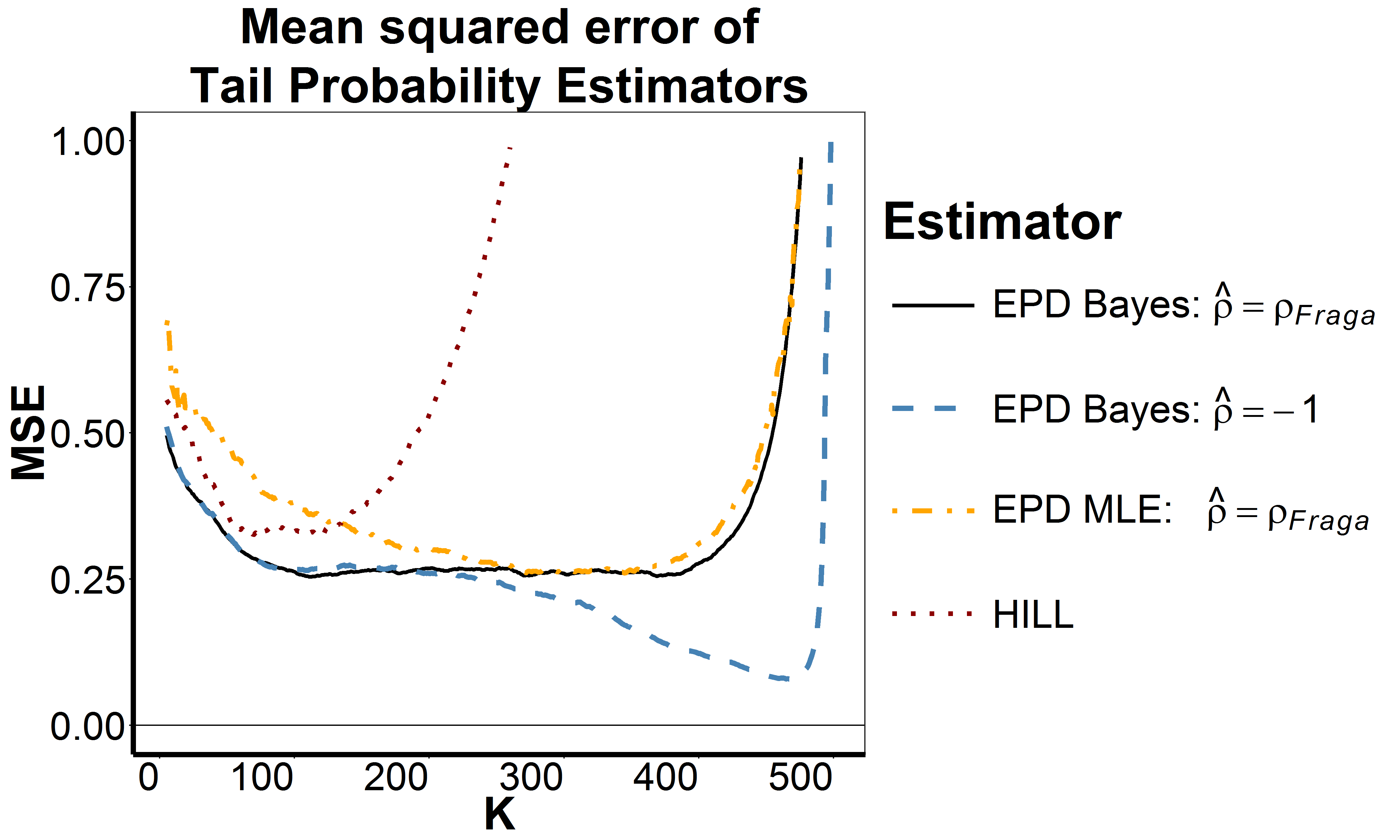}
 			\end{subfigure}	
 			\caption{Bias (top), variance (middle), and mean squared error (bottom) of the estimates of $\xi$ (left) and estimates of the exceedance probability (right) using the Bayesian estimates with $\hat\rho=\rho_{Fraga}$ and $\rho=-1$, the classical Hill-Weissman approach, and the EPD-ML approach, in the case of a \textbf{Fréchet} distribution with $\xi=0.5$ and $\rho=-1$.}
 			\label{c5}
 		\end{figure}
 		\begin{figure}[H]
 		\centering
 		\begin{subfigure}[h]{0.47\textwidth}
 		\includegraphics[width=80mm]{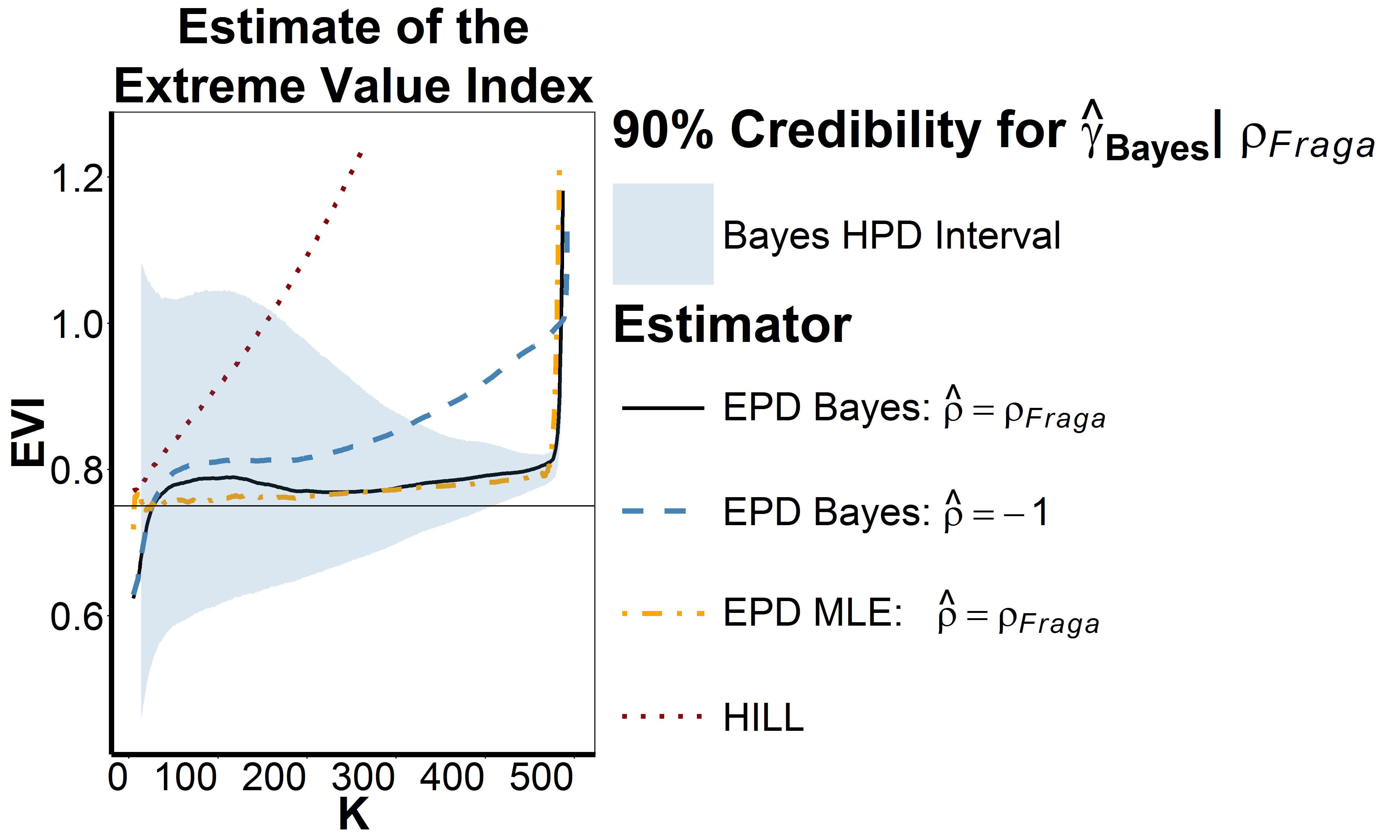}
 		\end{subfigure}
 		\hspace{\fill}
 		\begin{subfigure}[h]{0.45\textwidth}
 		\includegraphics[width=85mm]{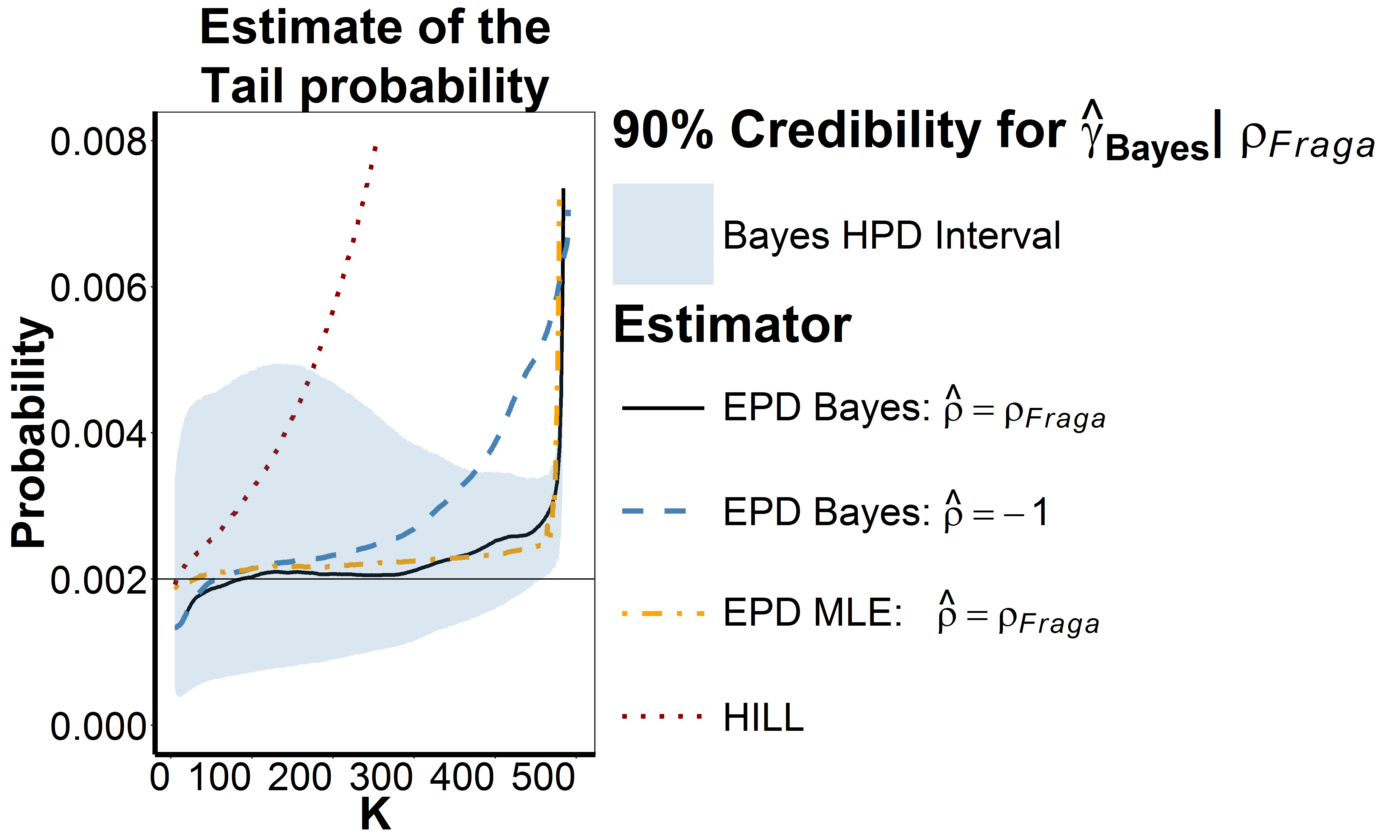}
 		\end{subfigure}
 		\bigskip
 		
 		\begin{subfigure}[h]{0.45\textwidth}
 		\includegraphics[width=75mm]{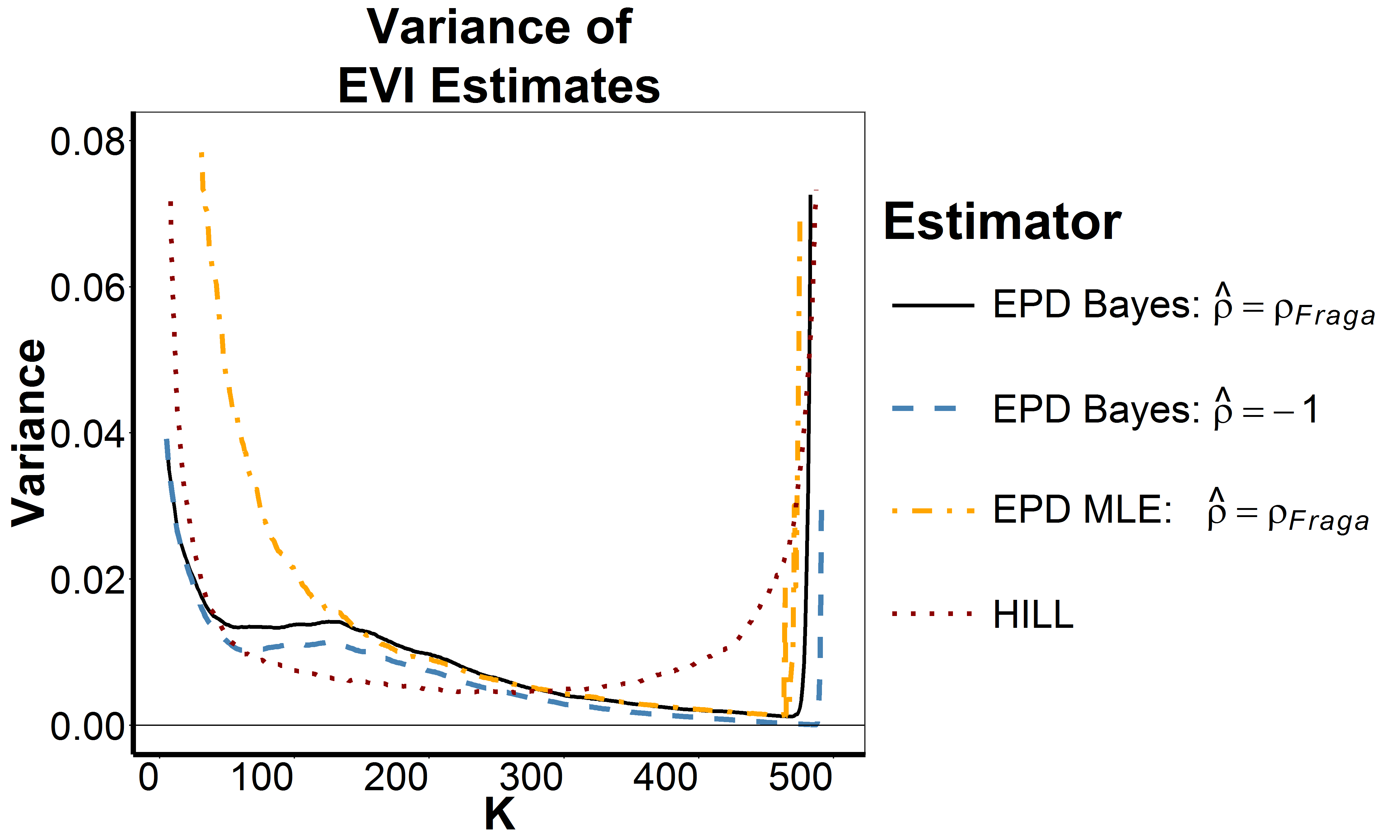}
 		\end{subfigure}	
 		\hspace{\fill}
 		\begin{subfigure}[h]{0.45\textwidth}
 		\includegraphics[width=75mm]{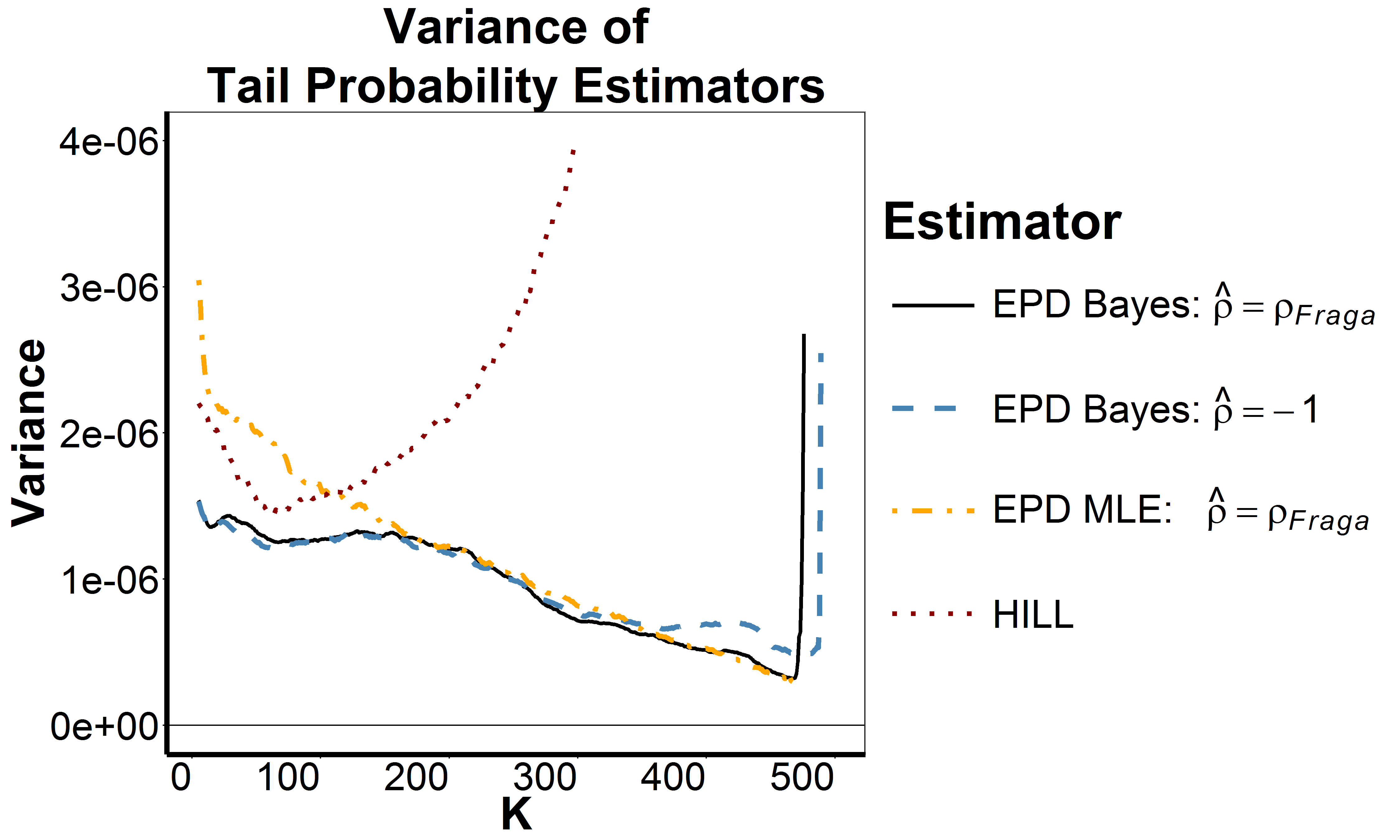}
 		\end{subfigure}	
 		  \bigskip

 		\begin{subfigure}[h]{0.45\textwidth}
 		\includegraphics[width=75mm]{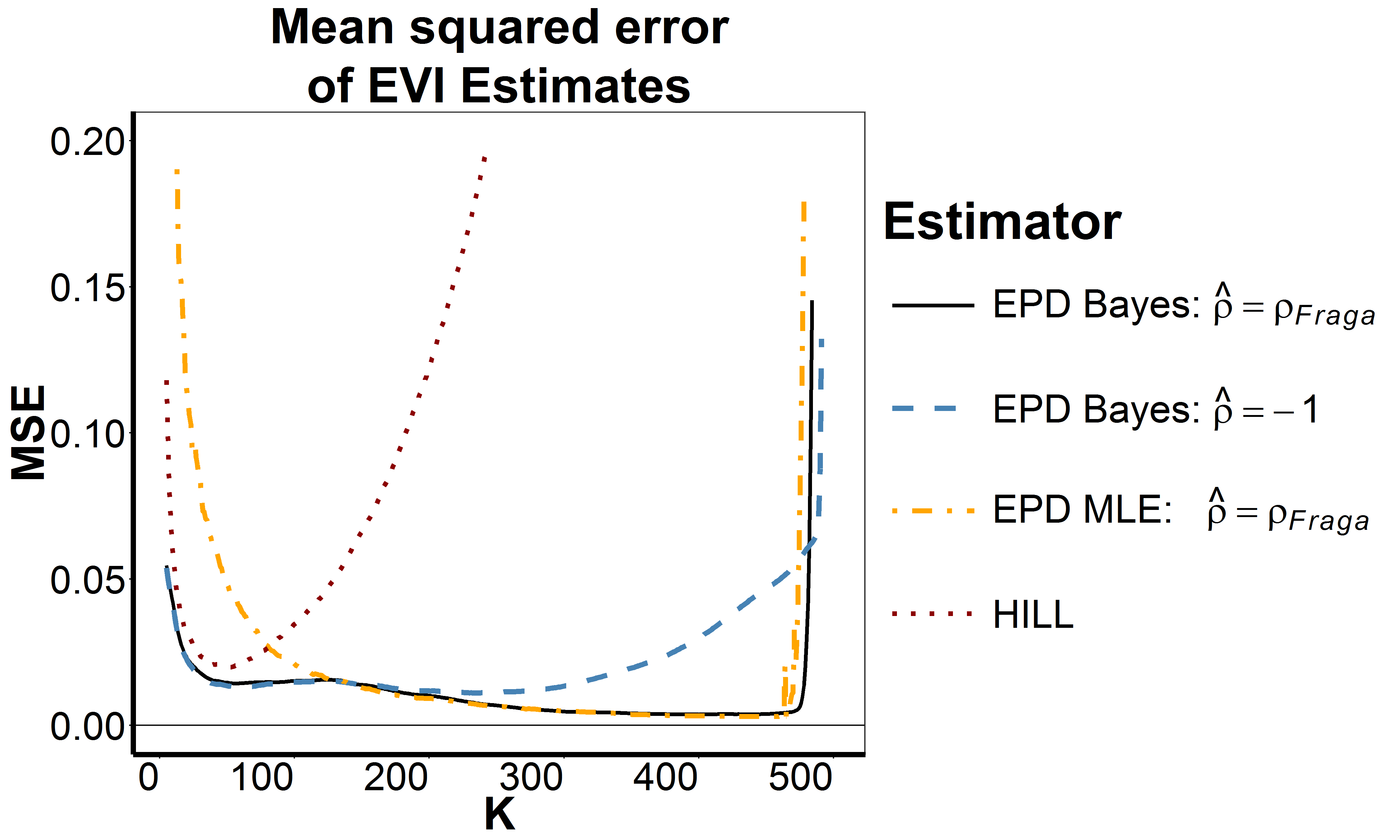}
 		\end{subfigure}	
 		\hspace{\fill}
 		\begin{subfigure}[h]{0.45\textwidth}
 		\includegraphics[width=75mm]{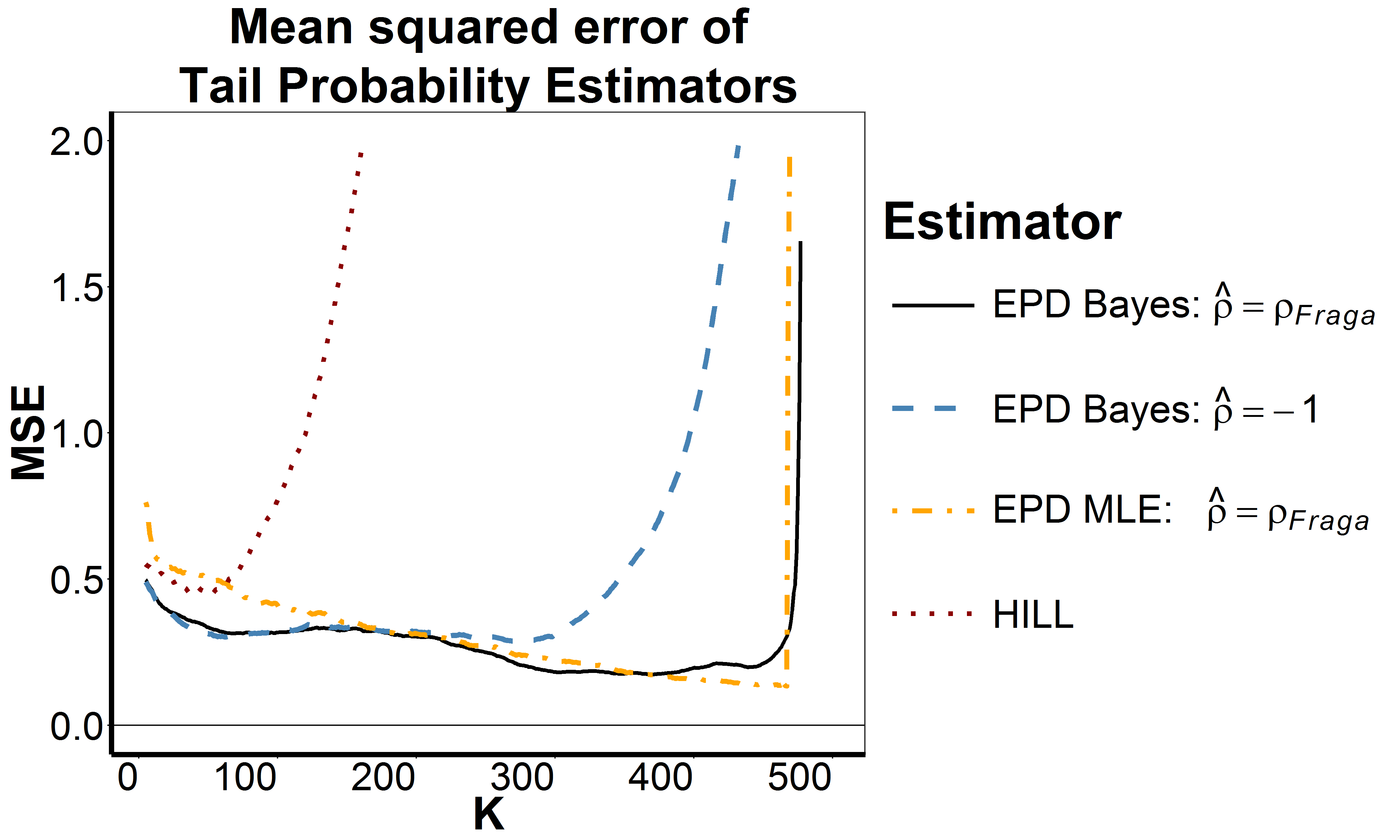}
 		\end{subfigure}	
 		\caption{Bias (top), variance (middle), and mean squared error (bottom) of the estimates of $\xi$ (left) and estimates of the exceedance probability (right) using the Bayesian estimates with $\hat\rho=\rho_{Fraga}$ and $\rho=-1$, the classical Hill-Weissman approach, and the EPD-ML approach, in the case of a \textbf{Burr} distribution with $\xi=0.75$ and $\rho=-0.75$.}
 		\label{c4}
 	\end{figure}
 		\begin{figure}[H]
		\centering
		\begin{subfigure}[h]{0.47\textwidth}
		\includegraphics[width=80mm]{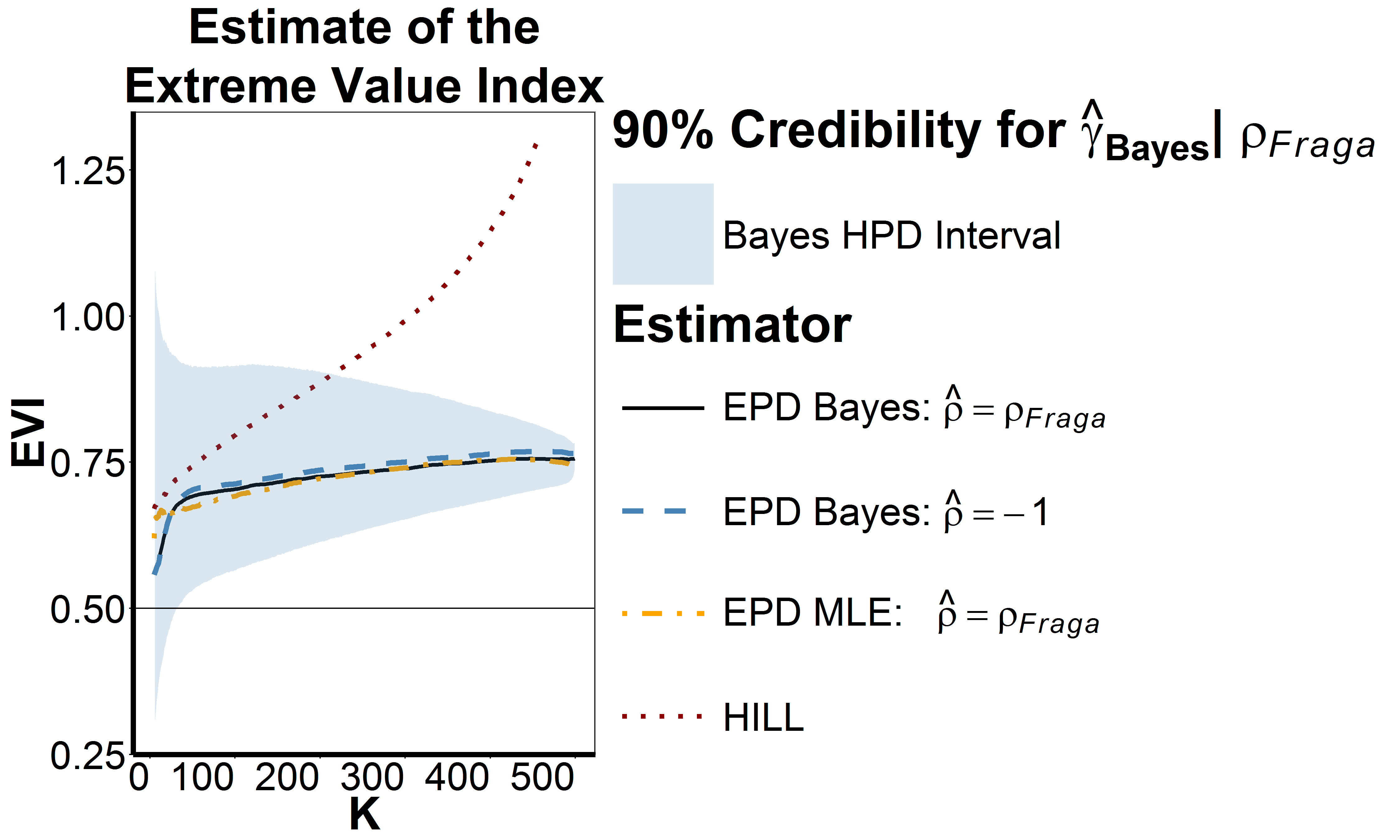}
		\end{subfigure}
		\hspace{\fill}
		\begin{subfigure}[h]{0.45\textwidth}
		\includegraphics[width=85mm]{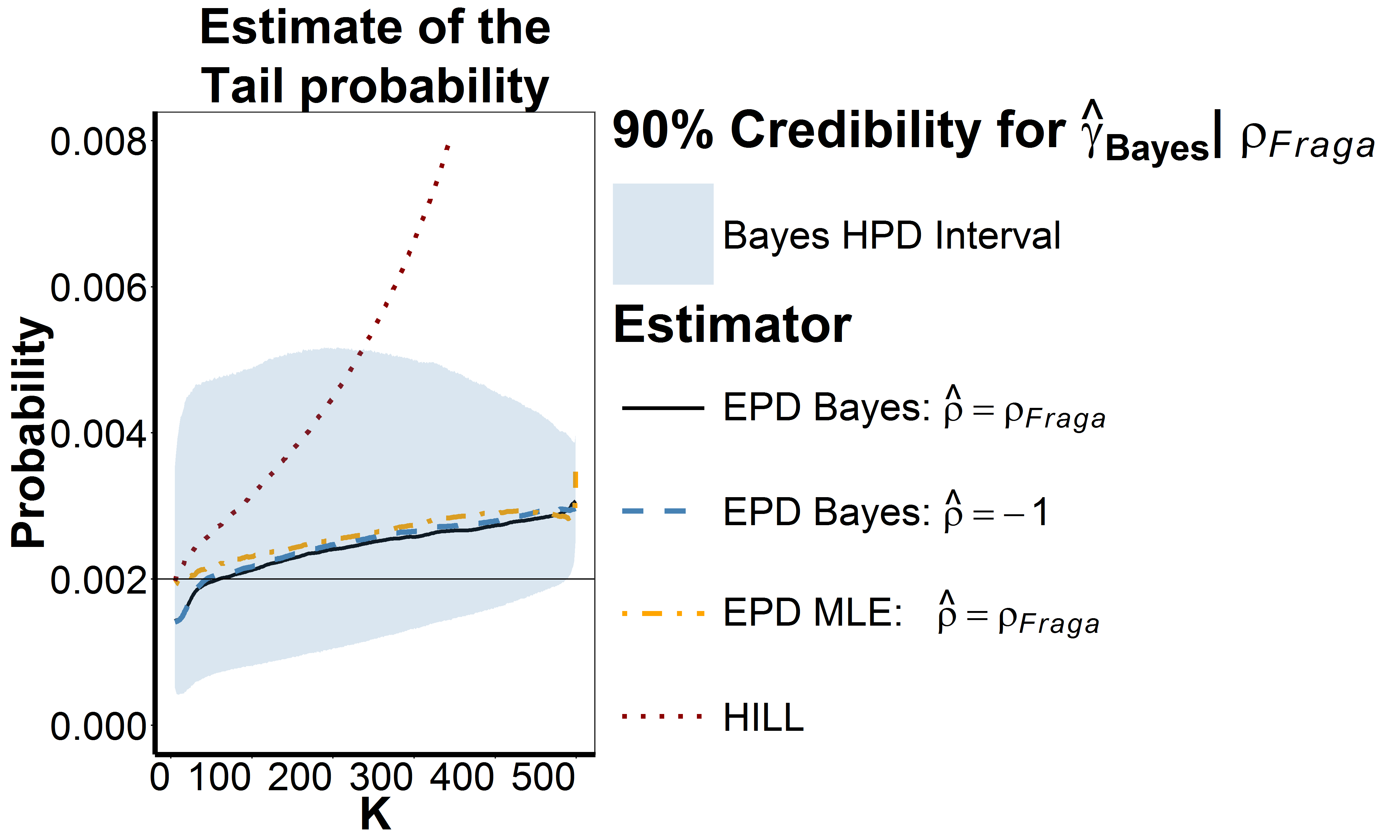}
		\end{subfigure}
		\bigskip
		
		\begin{subfigure}[h]{0.45\textwidth}
		\includegraphics[width=75mm]{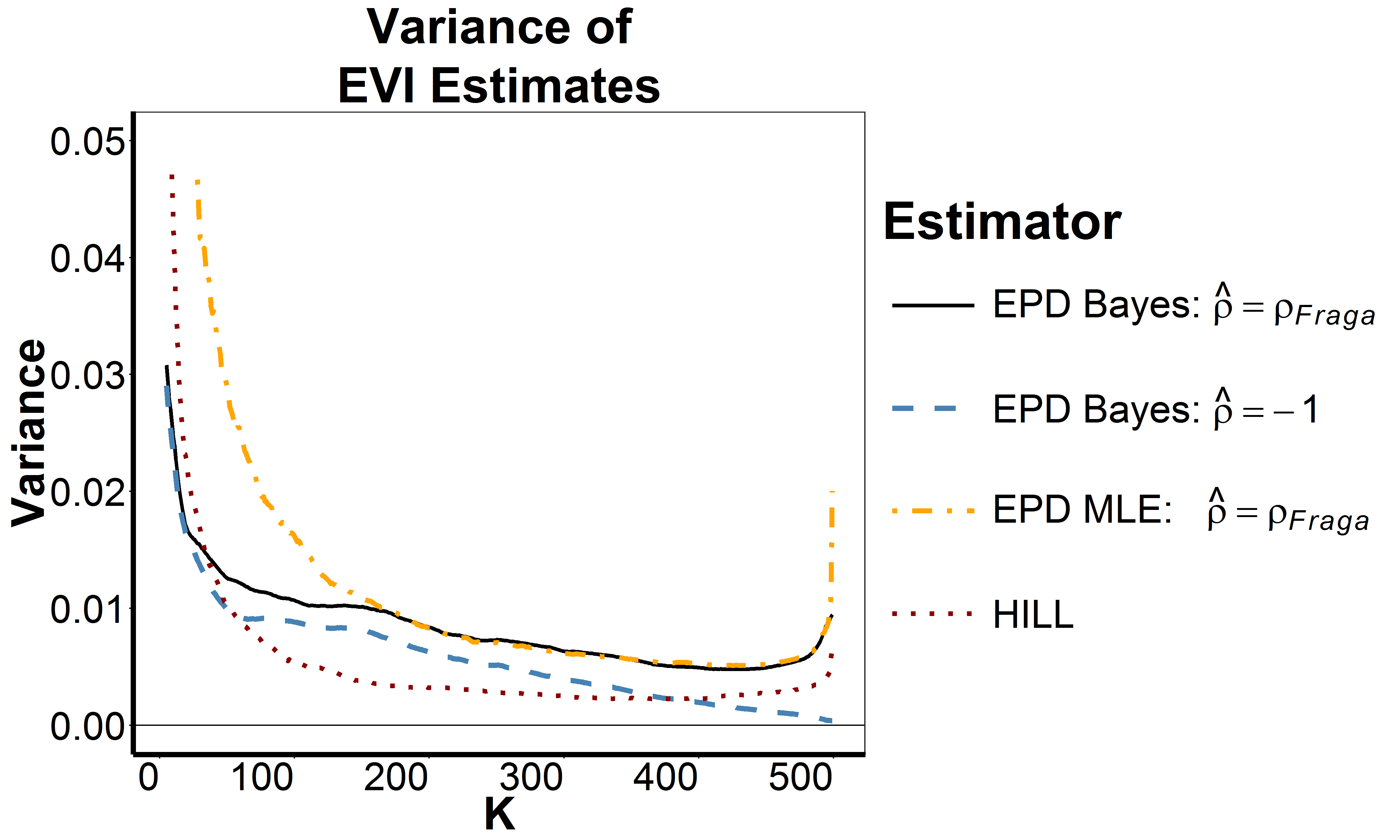}
		\end{subfigure}	
		\hspace{\fill}
		\begin{subfigure}[h]{0.45\textwidth}
		\includegraphics[width=75mm]{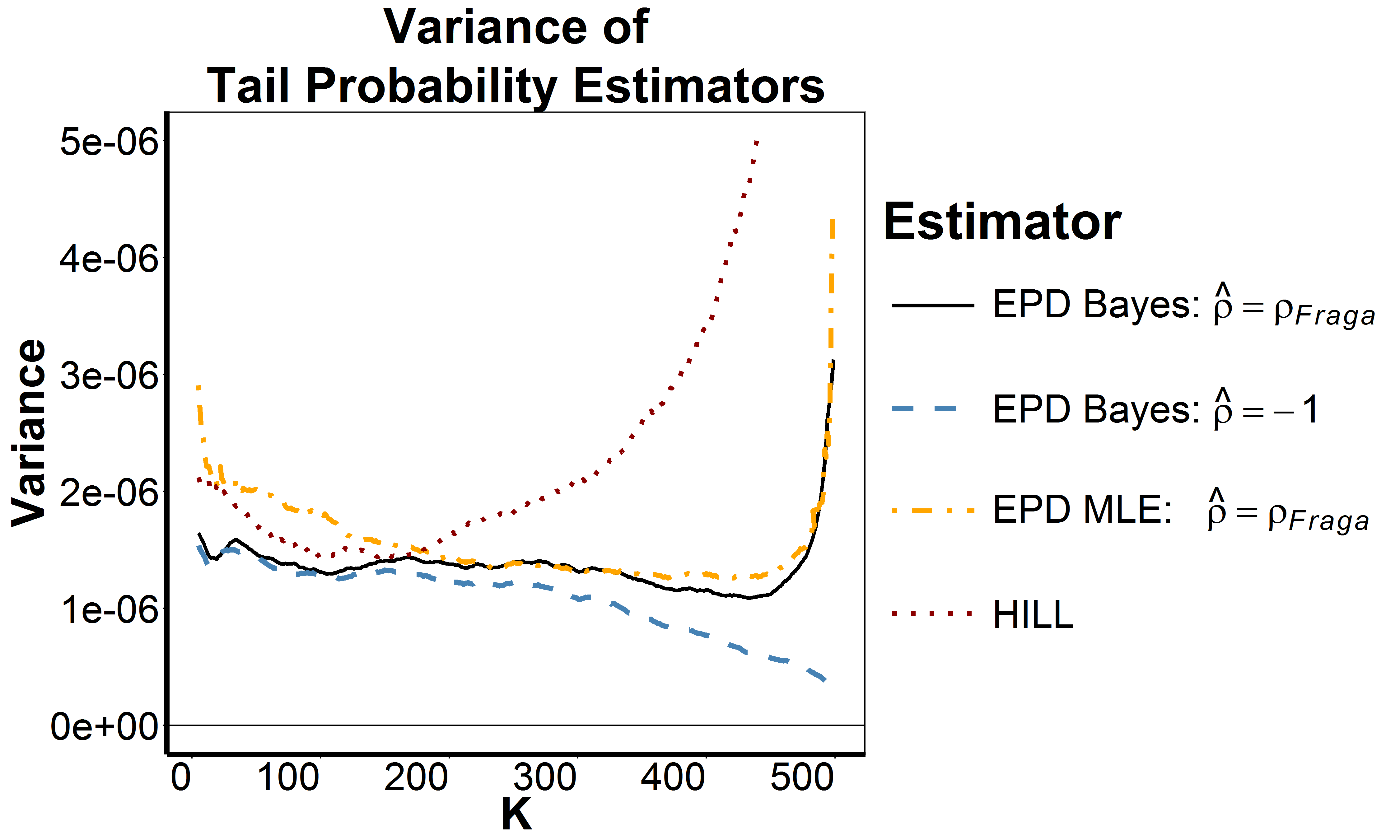}
		\end{subfigure}	
		\bigskip
		
		\begin{subfigure}[h]{0.45\textwidth}
		\includegraphics[width=75mm]{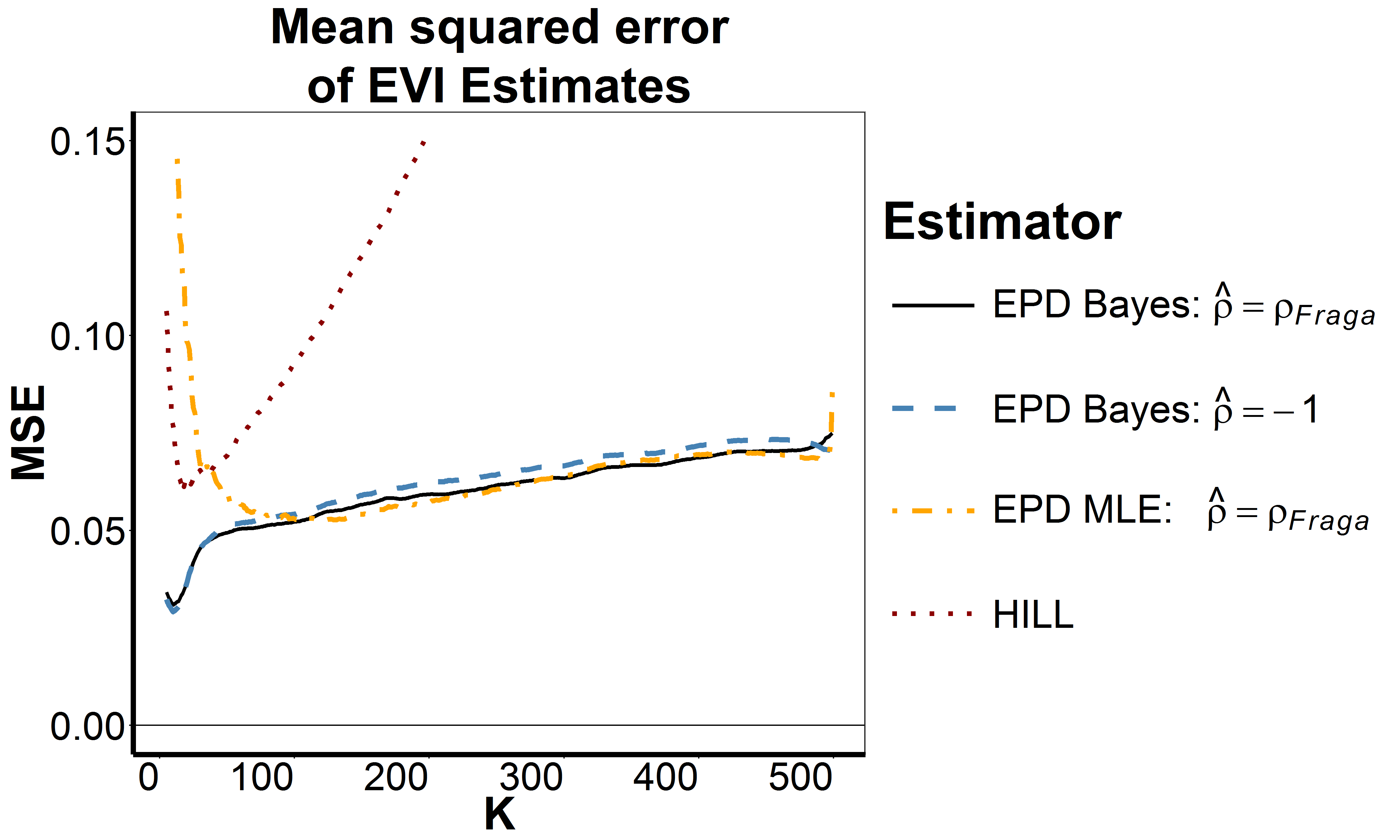}
		\end{subfigure}	
		\hspace{\fill}
		\begin{subfigure}[h]{0.45\textwidth}
		\includegraphics[width=75mm]{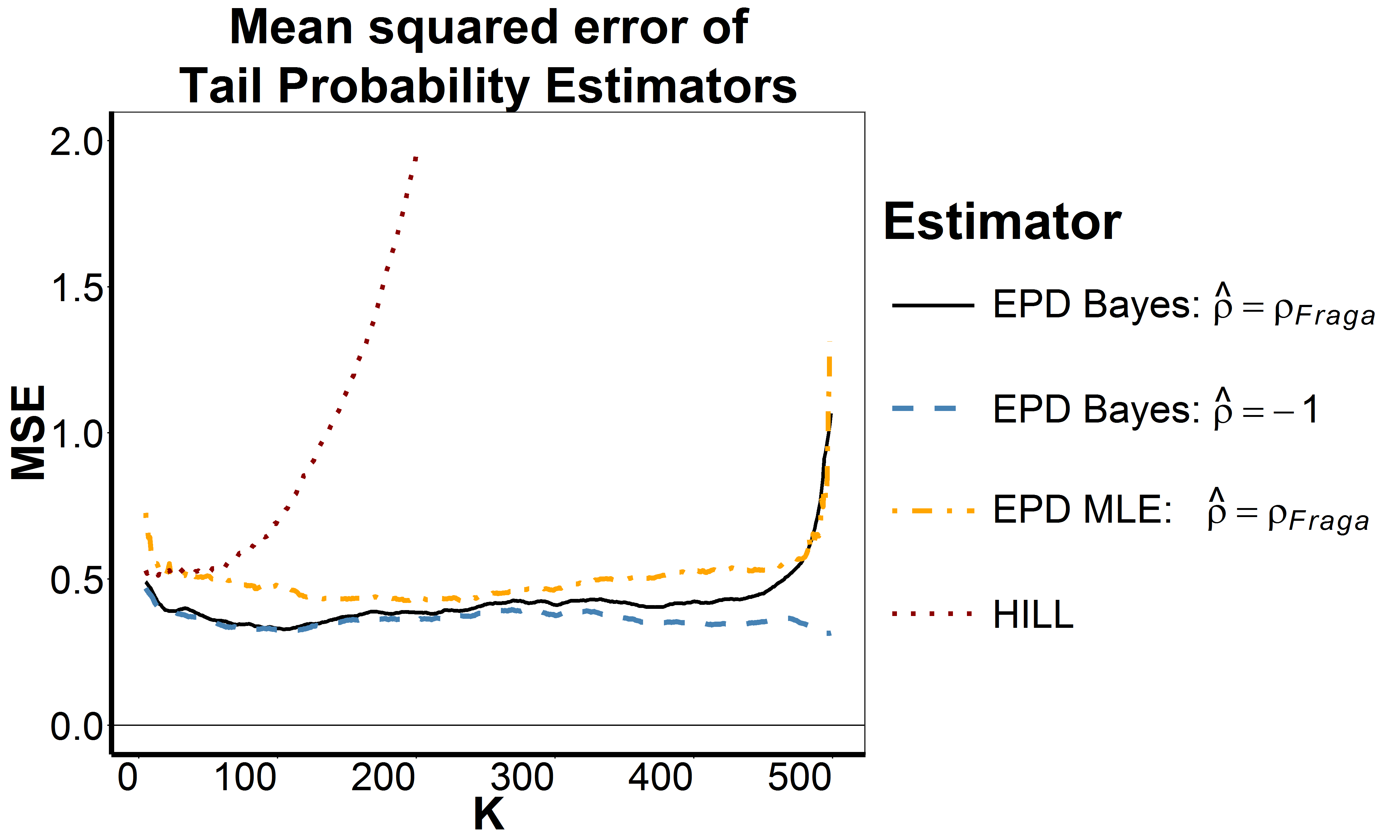}
		\end{subfigure}	
		\caption{Bias (top), variance (middle), and mean squared error (bottom) of the estimates of $\xi$ (left) and estimates of the exceedance probability (right) using the Bayesian estimates with $\hat\rho=\rho_{Fraga}$ and $\rho=-1$, the classical Hill-Weissman approach, and the EPD-ML approach, in the case of a \textbf{loggamma} distribution with $\xi=0.5$ and $\rho=0$.}
		\label{c6}
	\end{figure}
 
We conclude  that the finite sample behaviour of the proposed estimators follows the characteristics predicted by the asymptotic analysis to a great extent. For small $k$ the Bayesian estimators
$\xi_{k,n}^*$ and $\hat{p}^*_{x,k}$ show a similar behaviour as the Hill and Weissman estimators, while for larger $k$ the proposed estimators tend to follow the characteristics of the bias reduced EPD-ML estimator. The MSE of the Bayesian estimator of $P(X>x)$ is smallest, uniformly over the whole $k$ range and in all cases presented. Concerning the estimation of $\xi$, only in the Fréchet case and for small $k$ does the Hill estimator show a smaller MSE than the Bayesian estimator, while $\hat{\xi}^{*}_{k,n}$ then still shows a much smaller MSE than the EPD-ML estimator. \\ Also note that the version where the parameter $\rho$ is set to -1 does not differ too much from the use of the Fraga Alves {\it et al.} (2003) estimator. In case of the Fr\'echet distribution with $\rho=-1$, fixing this second order parameter at the correct value naturally yields some improvement in MSE, especially at large values of $k$. Similarly, fixing $\rho$ at an incorrect value, as in the case of the selected Burr distribution, yields larger MSE values  at large values of $k$. \\ Finally the results in case of the loggamma distribution are quite good. Hence it appears that the proposed method exhibits some robustness against deviations from the underlying model.
\\\\
In order to illustrate the use of the proposed method we consider the  weekly negative log-returns for Barclays PLC as studied in Reynkens {\it et al.} (2015). This data is divided in two sets, before and after the 2007 financial crisis:
\begin{itemize}
	\item Pre-Crisis: from January 1,1994 to August 7, 2007,
	\item Post-Crisis: from August 8,2007 to September 23, 2014.
\end{itemize}
We consider the estimation of $P(X>0.18)$ in both periods.  Whereas daily return data may suffer from serial dependence such as volatility clustering, such dependence is at least much weaker in weekly returns. We here use the proposed technique as a data analytic tool in order to assist the user in finding the  relevant level of the estimate as a function of $k$. Up to the start of the crisis an 18\% loss or higher appears to return on average once in 400 weeks (one-sided exceedance probability 0.005) based on $k/n \leq 0.25$, while as from the start of the crisis this return level decreases to close to once in 50 weeks (one-sided exceedance probability 0.04) based on $k/n \leq 0.50$. 
 \begin{figure}[H]
 		\centering
 		\begin{subfigure}[b]{1\textwidth}
 			\includegraphics[width=160mm]{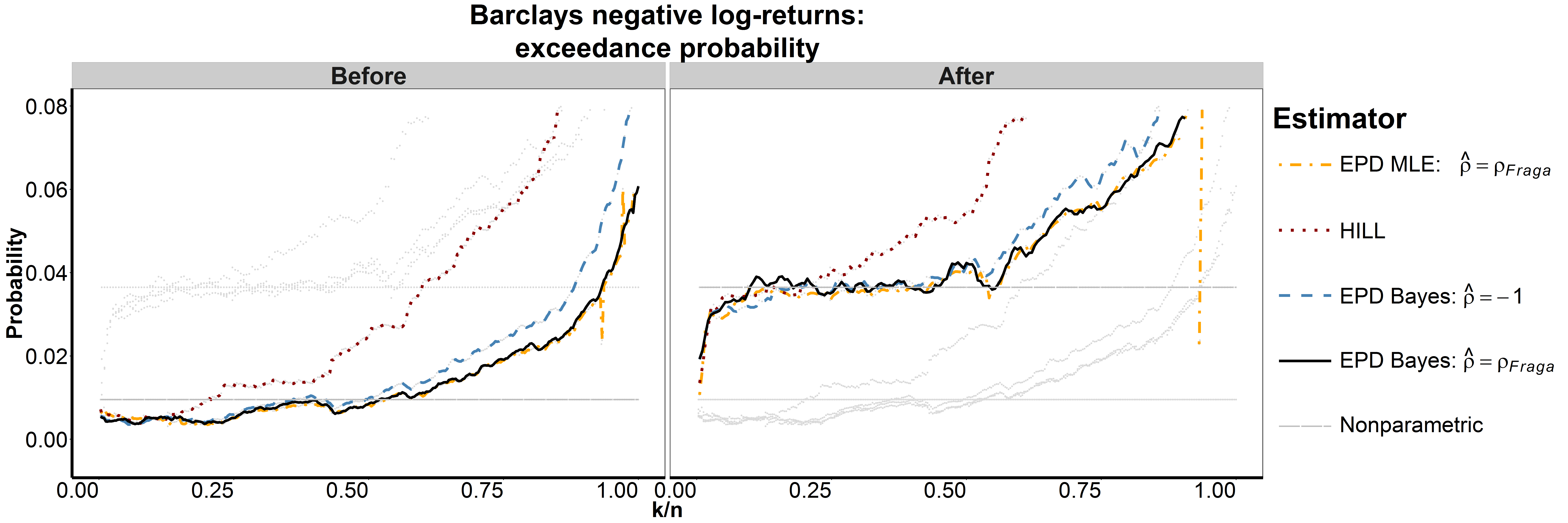}
 		\end{subfigure}
 		\caption{\footnotesize Estimates of  the  probability of exceeding 0.18 as a function of k/n for pre-crisis (left) and post-crisis (right) weekly negative log-returns of Barclays PLC}
 		\label{Barclays} 	
 \end{figure}

\section{Conclusion}
We proposed a simple adaptation to the bias reduction technique in tail estimation based on the EPD, which yields interesting MSE behaviour. In fact for larger thresholds the proposed estimators follow the behaviour of the classical Hill and Weissman estimators with small  bias and minimal variance, while the new estimators are never worse than the ML estimators based on the EPD approximation. In contrast to existing minimum variance bias reduced estimators which use third order tail conditions, the conditions can be kept minimal. In fact, setting the second order parameter $\hat \rho$ at $-1$ still yields acceptable results. 

In future work extensions of this approach  to other tail estimation problems will be investigated.
 
\noindent
{\bf Acknowledgment.} The authors take pleasure in thanking S. van der Merwe for his valuable advice concerning the Bayesian implementations.

\section{Appendix}
{\it Derivation of the expressions of ($\xib,\db$).} 
First consider the asymptotic approximations of the mode-posterior estimator of $\xi$ based on maximization of \eqref{logpost} with \eqref{pixi} and \eqref{pidelta}. From \eqref{logpost}-\eqref{pidelta} using expansions in $\delta \to 0$ 
we obtain 
\begin{eqnarray*}
{1 \over k}\log \pi(\xi,\delta |{\bf y}) &=& 
-(1+{1 \over k}) \log \xi -{1 \over k}(1+\xi) -(\frac{1}{\xi}+1){1 \over k}\sum_{j=1}^k\log y_{j,k} \\
&& -{\delta \over 1+\xi}{1 \over k}\sum_{j=1}^k (1-y_{j,k}^{\tau})
+ \delta  {1 \over k}\sum_{j=1}^k (1- (1+\tau)y_{j,k}^{\tau})\\
&&- {\delta^2 \over 2k\sigma_{k,n}^2}
+{\delta^2 \over 2(1+\xi)}{1 \over k}\sum_{j=1}^k (1-y_{j,k}^{\tau})^2
- {\delta^2 \over 2} {1 \over k}\sum_{j=1}^k (1- (1+\tau)y_{j,k}^{\tau})^2 \\
&&+ O(\delta^3) + c,
\end{eqnarray*}
where $c$ is a constant only depending on $\sigma_{k,n}^2$ and $\tau$. Also note that ${1 \over k}\sum_{j=1}^k\log y_{j,k}=H_{k,n}$.
 Then
 \begin{eqnarray*}
 {\partial \over \partial  \xi}{1 \over k}\log \pi(\xi,\delta |{\bf y}) &=& -{1 \over \xi} + {1 \over \xi^2} H_{k,n}
 +{\delta \over \xi^2}{1 \over k}\sum_{j=1}^k (1-y_{j,k}^{\tau})
 + O(\delta^2) + O({1 \over k}), \\
 {\partial \over \partial  \delta}{1 \over k}\log \pi(\xi,\delta |{\bf y}) &=& -{1 \over \xi}\left( 1- (1-\xi\tau){1 \over k}\sum_{j=1}^k y_{j,k}^{\tau} \right) -{\delta \over k\sigma^2_{k,n}}\\
 &&
 + {\delta \over \xi} \left( 1- 2(1-\xi\tau){1 \over k}\sum_{j=1}^k y_{j,k}^{\tau} 
 + (1-2\xi\tau-\xi\tau^2){1 \over k}\sum_{j=1}^k y_{j,k}^{2\tau}\right)+ O(\delta^2) + O({1 \over k}).
 \end{eqnarray*}
 
\vspace{0.4cm}\noindent 
{\it Derivation of Theorem.}
 Assuming $k\sigma_{k,n}^2 \to \mu >0$ as $k,n \to \infty$, $k/n \to 0$ and $\sqrt{k} a(n/k) \to \lambda$ we find using
$E_{k,n}(s) \to 1/(1-\xi s)$ (see Theorem A.1 in Beirlant et al., 2009)  that
  $$
    D^{\mathcal{B}}_{k,n} = -{\xi \over \mu} 
   + {\rho^4 \over \xi (1-2\rho)(1-\rho)^2} + o_p(1).
 $$
 Then, proceeding as in the proof of Theorem 3.1 in Beirlant et al. (2009), we obtain with $\Gamma_{k,n} = \sqrt{k}( H_{k,n}-\xi)$, $\mathbb{E}_{k,n}(s) =\sqrt{k}(E_{k,n}(s)-{1\over 1-\xi s })$ ($s<0$), that
 \begin{eqnarray*}
 \sqrt{k} \left( \xib-\xi \right)&=&
 \sqrt{k} \left( H_{k,n}-\xi - \db {\rho \over 1-\rho} \right) 
 \\
 &=& \Gamma_{k,n}-{\rho \over 1-\rho}\sqrt{k}\db \\
 &=& \Gamma_{k,n}\left(1+ {\rho^2 \over \xi(1-\rho^2)}
 {1 \over (\xi/\mu)+ \rho^4/\xi(1-2\rho)(1-\rho)^2} \right) \\
 &&- {\rho \over (\xi/\mu)+\rho^4/\xi(1-2\rho)(1-\rho)^2} \mathbb{E}_{k,n}(\hat\tau) +o_p(1)\\
 &=& \Gamma_{k,n} \left( 1+ {\rho^2 (1-2\rho)\over \zeta +\rho^4}\right) +  
 \mathbb{E}_{k,n}(\hat\tau)\left( {(-\rho)\xi(1-2\rho)(1-\rho)^2 \over \rho^4 + \zeta }\right)+o_p(1),
 \end{eqnarray*}
 where $\zeta= \xi^2(1-2\rho)(1-\rho)^2/\mu$. 
 Using Theorem A.1 in Beirlant et al. (2009), \eqref{Abias} and \eqref{Avar} follow under  $\sqrt{k} a(n/k) \to \lambda$.
\end{document}